\documentclass[12pt]{amsart}

\usepackage{amsmath,amsfonts,amscd,amsthm,amssymb,latexsym}

\font\teneufm=eufm10 \font\seveneufm=eufm7 \font\fiveeufm=eufm5
\newfam\frakturfam
\textfont\frakturfam=\teneufm \scriptfont\frakturfam=\seveneufm
\scriptscriptfont\frakturfam=\fiveeufm

\newtheorem{lm}{Lemma}
\newtheorem{theor}{Theorem}
\newtheorem{co}{Corollary}

\def\bee{\begin{eqnarray}}
\def\bes{\begin{eqnarray*}}
\def\eee{\end{eqnarray}}
\def\ees{\end{eqnarray*}}

\def\a{\alpha}

\def\Proof{{\sl Proof.}\ }

\title{The Freiheitssatz for Novikov algebras}

\begin{document}
\date{}
\maketitle

\begin{center}

{\bf Leonid Makar-Limanov}\footnote{Supported
by an NSA grant H98230-09-1-0008, by an NSF grant DMS-0904713, and a Fulbright fellowship awarded by the United States--Israel Educational Foundation; The Weizmann Institute of Science, Rehovot, Israel, University of Michigan, Ann Arbor, and
 Wayne State University, Detroit, MI 48202, USA,
e-mail: {\em lml@math.wayne.edu}}
and
{\bf Ualbai Umirbaev}\footnote{Supported by an NSF grant DMS-0904713 and by a grant of Kazakhstan; Eurasian National University,
 Astana, Kazakhstan and
 Wayne State University,
Detroit, MI 48202, USA,
e-mail: {\em umirbaev@math.wayne.edu}}

\end{center}

\begin{abstract}
We prove the Freiheitssatz for Novikov algebras in characteristic zero. It is also proved that
the variety of Novikov algebras is generated by a Novikov algebra on the space of polynomials $k[x]$ in a single variable $x$ over a field $k$ with respect to the multiplication $f\circ g=\partial(f)g$.
It follows that the base rank of the variety of Novikov algebras equals 1.
\end{abstract}


\noindent


\section{Introduction}

\hspace*{\parindent}

In 1930 W.\, Magnus proved
one of the most important theorems of the combinatorial group
theory (see \cite{Magnus}): \textsl{Let $G=\langle  x_1,x_2,\ldots,
x_n | r=1\rangle$ be a group defined by a single cyclically
reduced relator $r$. If $x_n$ appears in $r$, then the subgroup of
$G$ generated by $x_1,\ldots, x_{n-1}$ is a free group, freely
generated by $x_1,\ldots, x_{n-1}$.} He called it
\textit{the Freiheitssatz} (``freedom/independence theorem" in German). In the same paper W.\, Magnus proved the decidability of the word problem for
groups with a single defining relation. The Freiheitssatz for solvable and nilpotent
groups was researched by many authors (see, for example \cite{Romanovskii}).

In 1962 A.\,I.\,Shirshov \cite{Shir2} established the Freiheitssatz for Lie algebras
and proved the decidability of the word problem for Lie algebras
with a single defining relation. These results recently were generalized in \cite{KMLU}
for right-symmetric
algebras. In 1985 L.\,Makar-Limanov \cite{Makar2} proved the Freiheitssatz for
associative algebras of characteristic zero and in \cite{MLU2010} it was also proved for
Poisson algebras of characteristic zero. Note that the question of decidability of the word problem for
associative algebras and Poisson algebras with a single defining
relation and the Freiheitssatz for associative algebras in a
positive characteristic remain open. The Freiheitssatz for Poisson algebras in a
positive characteristic is not true \cite{MLU2010}.

In this paper we prove the Freiheitssatz for Novikov algebras over fields of characteristic zero.
There are two principal methods of proving the Freiheitssatz: one, employing
the combinatorics of free algebras, applied in \cite{Magnus,Romanovskii,Shir2,KMLU},
and the other, related to the study of algebraic and differential equations,
applied in \cite{Makar2,MLU2010}. The latter is used here.

Recall that an algebra $A$ over a field $k$ is called {\em
right-symmetric} if it satisfies the identity
\bee\label{f1}
(xy)z-x(yz)=(xz)y-x(zy).
\eee
In other words, the associator
$(x,y,z)=(xy)z-x(yz)$ is symmetric in $y$ and $z$. The variety of
right-symmetric algebras is Lie-admissible, i.e., each
right-symmetric algebra $A$ with the operation $[x,y]=xy-yx$ is a
Lie algebra.
 A right-symmetric algebra $A$ is called {\it Novikov}
 (\cite{Novikov}, \cite{Osborn}, \cite{Gelfand}),
 if it satisfies also the  identity
\bee\label{f2}
x(yz)=y(xz).
\eee

Let $k[x]$ be the polynomial algebra in a single variable $x$ over a field $k$ of characteristic $0$.
There are two interesting multiplications on   $k[x]$ (see, for example \cite{Dzhuma02,Dzhuma05,Dzhuma09}):
\bes
f*g=f\int_0^x g dx
\ees
and
\bes
f\circ g=\partial(f)g,  \ \ \partial=\frac{d}{d x}.
\ees
The algebra $\langle k[x],*\rangle$ is a free dual Leibniz algebra freely generated by $1$ and it was proved in \cite{NU} that
the variety of dual Leibniz algebras is generated by $\langle k[x],*\rangle$. The algebra $A=\langle k[x],\circ\rangle$ is a Novikov algebra \cite{Dzhuma02} and it is the main object of this paper. We prove that the variety of Novikov algebras is generated by $A$. It follows that the base rank of the variety of Novikov algebras is equal to $1$.

The paper is organized as follows. In Section 2 we prove that all identities of $A$ are corollaries of (\ref{f1})--(\ref{f2}).
In Section 3, using the homomorphisms of free Novikov algebras into $A$ and some results on differential equations from \cite{MLU2010}, we prove the Freiheitssatz.

\section{Identities}

\hspace*{\parindent}

Let $k$ be a field of characteristic $0$. Denote by $\mathfrak{N}$ the variety of Novikov algebras over $k$ and
denote by $\mathrm{N}\langle X\rangle$ the free Novikov algebra freely generated by $X=\{x_1,x_2,\ldots,x_n\}$. Put $x_1<x_2<\ldots<x_n$. In \cite{Dzhuma02,Dzhuma09} several constructions of a linear basis of $\mathrm{N}\langle X\rangle$ are given. We use a linear basis of $\mathrm{N}\langle X\rangle$ given in \cite{Dzhuma09} in
terms of Young diagrams.

Recall that a Young diagram is a set of boxes (we denote them by
bullets) with non-increasing numbers of boxes in each row. Rows
and columns  are numbered from the top to the bottom and from the left to the
right. Let $k$ be the  number of rows and $r_i$ be the number of
boxes in the $i$th row. The total number of boxes, $r_1+\cdots+r_k,$
is called the {\it degree} of the Young diagram.

To get a  Novikov diagram, we need to add
 one box (call it "a nose") to a Young diagram.
Namely, we need to add one more box to the first row, i.e.,
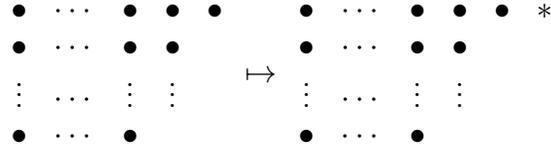
$$\begin{array}{ccccc} \bullet&\cdots&\bullet&\bullet&\bullet\\
\bullet&\cdots&\bullet&\bullet&\\ \vdots& \cdots&\vdots &\vdots&\\
\bullet&\cdots&\bullet&&\\
\end{array}
\mapsto
\begin{array}{cccccc}
\bullet&\cdots&\bullet&\bullet&\bullet&*\\
\bullet&\cdots&\bullet&\bullet&\\ \vdots&
\cdots&\vdots&\vdots&&\\\bullet&\cdots&\bullet&&&\\
\end{array}
$$ The number of boxes in a Novikov diagram is also called its {\it degree}. So, the
difference between the degrees of a Novikov diagram and the corresponding
Young diagram is  $1$.

 To construct
Novikov tableaux on $X$ we need to fill Novikov diagrams by
elements of $X$. Denote by $a_{i,j}$ the element of $X$
in the box $(i,j)$, that is the cross of the $i$-th row and the $j$-th
column. The {\it filling rules} are
\begin{itemize}
\item[(F1)] $a_{i,1}\ge a_{i+1,1},$ if $r_i=r_{i+1},
i=1,2,\ldots,k-1$;
\item[(F2)] the sequence of elements $a_{k,2},\ldots,a_{k,r_k},
a_{k-1,2},\ldots,a_{k-1,r_{k-1}},\ldots,a_{1,2},\ldots,
a_{1,r_1},a_{1,r_1+1}$ is non-decreasing.
\end{itemize}
In particular, all boxes beginning from the second place in each
row are labeled by non-decreasing elements of $X$.

 To any Novikov tableau
\bee\label{f3}
T=
\begin{array}{cccccc}
a_{1,1}&\cdots&\cdots&a_{1,r_1-1}&a_{1,r_1}&a_{1,r_1+1}\\
a_{2,1}&\cdots&a_{2,r_2-1}&a_{2,r_2}&\\ \vdots&
\cdots&\vdots&\vdots&&\\ a_{k,1}&\cdots&a_{k,r_k}&&&\\
\end{array}
\eee
associate a non-associative word
\bee\label{f4}
W_T= W_k(W_{k-1}(\ldots (W_2
W_1)\ldots )),
\eee
in the alphabet $X$ where
\bes
W_1=(\ldots
((a_{1,1} a_{1,2}) a_{1,3})\ldots  a_{1,r_1})
a_{1,r_1+1},
\ees
\bes
W_i=(\ldots
((a_{i,1} a_{i,2}) a_{i,3})\ldots
a_{i,r_i-1}) a_{i,r_i},\quad 1<i\leq k.
\ees

The set of all non-associative words associated with Novikov tableaux composes a linear basis of the free Novikov algebra $\mathrm{N}\langle X\rangle$ \cite{Dzhuma09}. 

Recall that $A=\langle k[x], \circ\rangle$ is the Novikov algebra on the space of the polynomial algebra $k[x]$ with respect to
multiplication $\circ$. For any $s=(s_1,\ldots,s_n)\in \mathbb{Z}_{+}^n$, where $\mathbb{Z}_{+}$ is the set of all nonnegative integers, we define a homomorphism
\bes
\overline{s} : \mathrm{N}\langle X\rangle \longrightarrow A=\langle k[x], \circ\rangle
\ees
given by $\overline{s}(x_i)=x^{s_i}$ for all $1\leq i\leq n$.

Consider the polynomial algebra $k[\lambda_1,\ldots,\lambda_n]$ in the variables $\lambda_1,\ldots,\lambda_n$. Put
 $\lambda=(\lambda_1,\ldots,\lambda_n)$ and $k[\lambda]=k[\lambda_1,\ldots,\lambda_n]$. Put also
 $x^{k[\lambda]}=\{x^{f(\lambda)} | f(\lambda)\in k[\lambda]\}$. Define a multiplication on $x^{k[\lambda]}$ by
 \bes
 x^{f(\lambda)}x^{g(\lambda)}=x^{f(\lambda)+g(\lambda)}.
 \ees
 Obviously, $x^{k[\lambda]}$ is a multiplicative copy of the additive group of $k[\lambda]$. Denote by $G$ the group algebra of  $x^{k[\lambda]}$ over $k[\lambda]$. It is easy to check that there exists a unique  $k[\lambda]$-linear derivation
 \bes
 D : G\longrightarrow G
 \ees
 such that $D(x^{f(\lambda)})=f(\lambda)x^{f(\lambda)-1}$ for all $f(\lambda)\in k[\lambda]$. With respect to
 \bes
 a\circ b= D(a)b, \ \ \ a,b\in G,
 \ees
 $G$ is a Novikov algebra again. Denote by $A(\lambda)$ the Novikov $k$-subalgebra of $G$ generated by
 $x^{\lambda_1},\ldots,x^{\lambda_n}$. The algebra $A(\lambda)$ looks like an algebra of general matrices (see, for example \cite{Drensky00}).

 Let
\bes
\overline{\lambda} : \mathrm{N}\langle X\rangle \longrightarrow A(\lambda)
\ees
 be an epimorphism of Novikov algebras defined by $\overline{\lambda}(x_i)=x^{\lambda_i}$ for all $1\leq i\leq n$. Note that $\overline{\lambda}$ is a "general" element for the set of all homomorphisms $\overline{s}$, where $\overline{s}\in \mathbb{Z}_{+}^n$. A homomorphism $\overline{s}$ is called a {\em specialization} of $\overline{\lambda}$.

Now we fix a Novikov tableau $T$ and its associated non-associative
word $W_T$ from (\ref{f3})--(\ref{f4}). Denote by $\deg$ the standard degree function on  $\mathrm{N}\langle X\rangle$ and by $\deg_{x_i}$ the degree function with respect to $x_i$ for all $1\leq i\leq n$. Denote by $d$ the degree of $T$ and by $d_i$ the number of occurrences of $x_i$ in $T$. Obviously, $d=\deg\,W_T$, $d_i=\deg_{x_i}W_T$, and
\bes
\overline{\lambda}(W_T)=f_T(\lambda) x^{g_T(\lambda)}
\ees
for some $f_T(\lambda), g_T(\lambda)\in k[\lambda]=k[\lambda_1,\ldots,\lambda_n]$.

Our first aim is to calculate the polynomials $f_T(\lambda)$ and $g_T(\lambda)$. For this reason we change the tableau $T$ from (\ref{f3}) by substituting $\lambda_i$ instead of $x_i$ for all $1\leq i\leq n$. Denote the new tableau by $T(\lambda)$. Then denote by $\lambda_{i,j}$ the element
in the box $(i,j)$ of $T(\lambda)$. In fact, we have just changed all $a_{i,j}$ to $\lambda_{i,j}$ in (\ref{f3}).
\begin{lm}\label{l1} The following statements are true:
\begin{itemize}
\item[(a)] $g_T(\lambda)=(d_1\lambda_1+\ldots+d_n\lambda_n-d+1)$;
\item[(b)] $f_T(\lambda)=f_1f_2\ldots f_k$
where
\bes
f_i=\lambda_{i,1}(\lambda_{i,1}+\lambda_{i,2}-1)\ldots(\lambda_{i,1}+\ldots+\lambda_{i,r_i}-r_i+1), \ 1\leq i\leq k.
\ees
\end{itemize}
\end{lm}
\Proof
Direct calculation gives that
\bes
\overline{\lambda}(W_1)=\overline{\lambda}((\cdots
((a_{1,1} a_{1,2}) a_{1,3})\cdots  a_{1,r_1})
a_{1,r_1+1})\\
=\overline{\lambda}((\cdots
((x^{\lambda_{1,1}}\circ x^{\lambda_{1,2}})\circ x^{\lambda_{1,3}})\circ\cdots \circ x^{\lambda_{1,r_1}})\circ
x^{\lambda_{1,r_1+1}})\\
=\lambda_{1,1}(\lambda_{1,1}+\lambda_{1,2}-1)\ldots(\lambda_{1,1}+\ldots+\lambda_{1,r_1}-r_1+1)
x^{(\lambda_{1,1}+\ldots+\lambda_{1,r_1}+\lambda_{1,r_1+1}-r_1)}.
\ees
Using this and leading an induction on $k$ we get
\bes
\overline{\lambda}(W_k)=\lambda_{k,1}(\lambda_{k,1}+\lambda_{k,2}-1)\ldots(\lambda_{k,1}+\ldots+\lambda_{k,r_k-1}-r_k+2)
x^{(\lambda_{k,1}+\ldots+\lambda_{k,r_k}-r_k+1)}
\ees
and
\bes
\overline{\lambda}(W_{k-1}(W_{k-2}\cdots (W_2
W_1)\cdots ))=f_1f_2\ldots f_{k-1}x^s,
\ees
where
$s=\sum_{i<k,j} \lambda_{i,j}-d+r_k+1$. Consequently,
\bes
\overline{\lambda}(W_T)=\overline{\lambda}(W_k)\circ \overline{\lambda}(W_{k-1}(W_{k-2}\cdots (W_2
W_1)\cdots ))\\
=\partial(\overline{\lambda}(W_k)) \overline{\lambda}(W_{k-1}(W_{k-2}\cdots (W_2
W_1)\cdots ))\\
=f_k x^{(\lambda_{k,1}+\ldots+\lambda_{k,r_k}-r_k)}f_1f_2\ldots f_{k-1}x^s=f_T x^t,
\ees
where $t=\lambda_{k,1}+\ldots+\lambda_{k,r_k}-r_k+s=\sum_{i,j} \lambda_{i,j}-d+1=g_T(\lambda)$.
$\Box$

\begin{lm}\label{l2}
A Novikov tableau $T$ is uniquely defined by the polynomials $f_T(\lambda)$ and $g_T(\lambda)$.
\end{lm}
\Proof For any linear form $l$ of the type
\bee\label{f5}
l=t_1\lambda_1+\ldots+t_n\lambda_n-t_1-\ldots-t_n+1
\eee
 we put
$\alpha(l)=t_1+\ldots+t_n$ and $\widehat{l}=t_1\lambda_1+\ldots +t_n\lambda_n$. Let $s_i$ be the number of boxes in the $i$-th column of the Young diagram corresponding to $T$. It follows from Lemma \ref{l1}(b) that $s_i$ is equal to the number of all divisors $l$ of $f_T$ of the form (\ref{f5}) with $\a(l)=i$, counted together with multiplicity. So, the Young diagram and the Novikov diagram corresponding to $T$ are uniquely defined.

By Lemma \ref{l1}(a), the degree of $T$ and the number of occurrences of $x_i$ in $T$ are also uniquely defined by $g_T(\lambda)$. It follows from Lemma \ref{l1}(b) that $x_i$ occurs in the first column of $T$ $m$-times if and only if $\lambda_i^m | f_T$ and $\lambda_i^{m+1} \nmid f_T$. Consequently, the elements of all columns of $T$, except the first one, are uniquely defined by the filling rule (F2).

So, the only question to answer is that how to arrange the elements of the first row. Let $l_1,\ldots,l_s$ be all divisors of $f_T$ of the form (\ref{f5}) with maximal $\alpha=\alpha(l_1)=\ldots=\alpha(l_s)$. By Lemma \ref{l1}(b),  $l_1,\ldots,l_s$ correspond to the first $s$ rows of $T$ and the first $s$ rows of the Young diagram corresponding to $T$ have lengths $r_1=\ldots=r_s=\alpha$. We have
\bes
\sum_{1\leq i\leq s} \sum_{1\leq j\leq r_i} \lambda_{i,j}=\widehat{l_1}+\ldots+\widehat{l_s}.
\ees
Suppose that
\bes
\sum_{1\leq i\leq s}  \lambda_{i,1}=\widehat{l_1}+\ldots+\widehat{l_s}-\sum_{1\leq i\leq s} \sum_{2\leq j\leq r_i} \lambda_{i,j}=\sum_{i=1}^n t_i\lambda_i.
\ees
Obviously $t_i\geq 0$, $t_1+\ldots+t_n=s$, and
\bes
(a_{1,1},\ldots,a_{s,1})=(\underbrace{x_n,\ldots,x_n}_{t_n},\ldots,\underbrace{x_1,\ldots,x_1}_{t_1})
\ees
by the filling rule (F1). So, the first $s$ rows of the Novikov tableaux $T$ are uniquely determined. Consequently,
the polynomials $f_1,\ldots,f_s$ are also uniquely determined. Using the polynomial $f_T/(f_1\ldots f_s)$ and continuing the same discussions, we can uniquely determine $T$.
$\Box$

Denote by $\mathbb{T}_n$ the set of all Novikov tableaux of degree $n$ on $X=\{x_1,\ldots,x_n\}$ without repeated elements. Then $\{W_T | T\in \mathbb{T}_n\}$ is a linear basis of the space of all multi-linear homogeneous of degree $n$ elements of the free Novikov algebra $\mathrm{N}\langle X\rangle$ \cite{Dzhuma09}.
\begin{co}\label{c1} Suppose that $T\in \mathbb{T}_n$. Then $T$ is uniquely defined by $f_T$.
\end{co}

Let $u=\lambda_1^{k_1}\ldots \lambda_n^{k_n}$ be an arbitrary monomial in $k[\lambda]=k[\lambda_1,\ldots,\lambda_n]$. Put $|u|=k_1+\ldots+k_n$. Put also $\gamma(u)=(s_1,\ldots,s_n)$ if $u=\lambda_{\sigma(1)}^{s_1}\ldots \lambda_{\sigma(n)}^{s_n}$ where $\sigma$ is a permutation on $\{1,\ldots,n\}$ and $s_1\geq s_2\geq \ldots\geq s_n$. We define a linear order $\preceq$ on the set of all monomials of $k[\lambda]$. If $u$ and $v$ are two monomials then put $u\preceq v$ if $|u|<|v|$ or  $|u|=|v|$ and $\gamma(u)$ is preceeds to $\gamma(v)$ with respect to the lexicographical order (from left to right) on $\mathbb{Z}_+^n$. If $|u|=|v|$ and $\gamma(u)=\gamma(v)$ then
$u\preceq v$ is defined arbitrarily.
 For any $f\in k[\lambda]$ denote by $\widetilde{f}$ its highest term with respect to $\preceq$.

 The statement of the next corollary trivially follows from Lemma \ref{l1}(b).
\begin{co}\label{c2} Suppose that $T\in \mathbb{T}_n$ and $(a_{1,1},a_{2,1},\ldots,a_{k,1})=(x_{i_1},x_{i_2},\ldots,x_{i_k})$ in (\ref{f3}).  Then,
\bes
\widetilde{f_T}= \lambda_{i_1}^{r_1}\lambda_{i_2}^{r_2}\ldots \lambda_{i_k}^{r_k} \ \ {\it and} \ \
\gamma(\widetilde{f_T})= (r_1,r_2,\ldots,r_k).
\ees
\end{co}

\begin{co}\label{c3} The set of polynomials $f_T\in k[\lambda]$, where $T$ runs over $\mathbb{T}_n$, is linearly independent over $k$.
\end{co}
\Proof
 Suppose that $(a_{1,1},a_{2,1},\ldots,a_{k,1})=(x_{i_1},x_{i_2},\ldots,x_{i_k})$ in (\ref{f3}). Then,
  $\gamma(\widetilde{f_T})= (r_1,r_2,\ldots,r_k)$ by Corollary \ref{c2}.   It follows that the Novikov diagram corresponding to $T$ is uniquely determined by $\widetilde{f_T}$. Moreover, $x_{i_s}$ is the first element of the row with length $r_s$. Then the filling rule (F1) uniquely determines the elements of the first row of $T$. The filling rule (F2) determines uniquely the other part of $T$.

So, the mapping $T\mapsto \widetilde{f_T}$ associates different tableaux to different basis elements of $k[\lambda]$. Consequently, the set of polynomials $\widetilde{f_T}$, where $T$ runs over $\mathbb{T}_n$, is linearly independent. This proves the lemma. $\Box$

In characteristic $0$ any identity is equivalent to the set of multi-linear homogeneous identities \cite{KBKA}.
Any nontrivial multi-linear homogeneous  Novikov identity of degree $n$ can be written as
\bee\label{f6}
\sum_{T\in \mathbb{T}_n} \a_T W_T=0
\eee
where $\a_T\in k$ and at least one of $\a_T$ is nonzero.

\begin{theor}\label{t1}
The Novikov algebra $A=\langle k[x],\circ\rangle$ does not satisfy any nontrivial Novikov identity.
\end{theor}
\Proof Suppose that  $A$ satisfies a nontrivial identity of the form (\ref{f6}). Consider the  homomorphism $\overline{\lambda}$. Applying $\overline{\lambda}$ to the left hand side of (\ref{f6}) we get
\bes
\overline{\lambda}(\sum_{T\in \mathbb{T}_n} \a_T W_T)=\sum_{T\in \mathbb{T}_n} \a_T f_T x^{g_T}
=(\sum_{T\in \mathbb{T}_n} \a_T f_T)x^{\lambda_1+\ldots+\lambda_n-n+1}
\ees
since $g_T(\lambda)=\lambda_1+\ldots+\lambda_n-n+1$ for all $T$.
By Corollary \ref{c3}, $\sum_{T} \a_T f_T$ is a nontrivial polynomial from $k[\lambda]$.
Then it is not difficult to find $s=(s_1,\ldots,s_n)\in \mathbb{Z}_{+}^n$ such that $\sum_{T} \a_T f_T(s_1,\ldots,s_n)\neq 0$. This means that the image of the left hand side of (\ref{f6}) under the homomorphism $\overline{s}$ is not equal to $0$.  Consequently, (\ref{f6}) is not a nontrivial identity of $A$. $\Box$
\begin{co}\label{c4} The variety of Novikov algebras $\mathfrak{N}$ is generated by $A=\langle k[x],\circ\rangle$, i.e.,
$\mathfrak{N}=\mathrm{Var}\,A$.
\end{co}

Recall that the least natural number $n$ such that the variety $Var(\mathrm{N}\langle x_1,x_2,\ldots,x_n\rangle)$ of algebras generated by  $\mathrm{N}\langle x_1,x_2,\ldots,x_n\rangle$ is equal to $\mathfrak{N}$ is called the {\em base rank}  $rb(\mathfrak{N})$ of the variety $\mathfrak{N}$ (see, for example \cite{NU}).

\begin{co}\label{c5}
The base rank of the variety of Novikov algebras is equal to one.
\end{co}
\Proof Consider the ideal $I$ of the polynomial algebra $k[x]$ generated by $x^2$. It is easy to check that $\langle I,\circ\rangle$ is a Novikov algebra generated by $x^2$. In the proof of Theorem \ref{t1}, we can easily chose $s=(s_1,\ldots,s_n)$ such that $s_i\geq 2$ for all $i$. Consequently, $\langle I,\circ\rangle$ does not satisfy any nontrivial Novikov identity. Then, $\mathfrak{N}=\mathrm{Var}\,\langle I,\circ\rangle$. We have   $\mathrm{Var}(\mathrm{N}\langle x_1\rangle)\supseteq \mathrm{Var}\,\langle I,\circ\rangle$  since $\langle I,\circ\rangle$ is a homomorphic image of $\mathrm{N}\langle x_1\rangle$. Therefore, $\mathfrak{N}=\mathrm{Var}(\mathrm{N}\langle x_1\rangle)$. $\Box$

\section{The Freiheitssatz}

\hspace*{\parindent}
To prove the Freiheitssatz we need the following corollary of Proposition 1 from \cite{MLU2010}.
\begin{co}\label{c6}\cite{MLU2010}
Let $f(x,t_{\a_1},t_{\a_2},\ldots,t_{\a_m})\in k[x,t_{\a_1},t_{\a_2},\ldots,t_{\a_m}]$ and $\a_1< \a_2< \ldots < \a_m$ are nonnegative integers.
Suppose that there exists $(c,c_{\a_1},c_{\a_2},\ldots,c_{\a_m})\in k^{1+m}$ so that
$f(c,c_{\a_1},c_{\a_2},\ldots,c_{\a_m})=0$ and $\frac{\partial f}{\partial t_{\a_m}}(c,c_{\a_1},c_{\a_2},\ldots,c_{\a_m})\neq 0$. Then the differential equation
\bes
f(x,\partial^{\a_1}(T),\partial^{\a_2}(T),\ldots,\partial^{\a_m}(T))=0
\ees
has a solution in the formal power series algebra $k[[x-c]]$.
\end{co}
Note that in the formulation of this corollary, the variables $x,t_{\a_1},t_{\a_2},\ldots,t_{\a_m}$ are independent variables, $\partial$ is the standard derivation $\frac{d}{d x}$ of $k[[x-c]]\supseteq k[x]$, and $\partial^{\a_i}$ is the $\a_i$th power of $\partial$.

If $f\in \mathrm{N}\langle x_1,\ldots,x_n\rangle$,  then we denote $\mathrm{id}(f)$ the ideal of $\mathrm{N}\langle x_1,\ldots,x_n\rangle$ generated by $f$.

\begin{theor}\label{t2}{\bf (Freiheitssatz)} Let $\mathrm{N}\langle x_1,\ldots,x_n\rangle$ be the free Novikov algebra over a field
$k$ of characteristic $0$ in the variables $x_1,\ldots,x_n$. If $f \in \mathrm{N}\langle x_1,\ldots,x_n\rangle$ and $f\notin
\mathrm{N}\langle x_1,\ldots,x_{n-1}\rangle$, then $\mathrm{id}(f)\cap \mathrm{N}\langle x_1,\ldots,x_{n-1}\rangle=0$.
\end{theor}
\Proof Without loss of generality we may assume that $k$ is algebraically closed and that
$f(x_1,\ldots,x_{n-1}, 0) \neq 0$.
The theorem will be proved if for $f$ and any nonzero $g\in \mathrm{N}\langle x_1,\ldots,x_{n-1}\rangle$ there exist a Novikov  algebra $B$ and
a homomorphism $\theta : \mathrm{N}\langle x_1,\ldots,x_n\rangle\rightarrow B$ of Novikov algebras such that
$\theta(g)\neq 0, \theta(f)=0$.

Let $\hat{f}$ be the highest homogeneous part of $f$ with respect to $x_n$.
By Theorem \ref{t1}, there exists a homomorphism
$\phi : \mathrm{N}\langle x_1,\ldots,x_n\rangle\rightarrow A=\langle k[x],\circ\rangle$ such that $\phi((gf)\hat{f})\neq 0$. Denote by $Z_1,Z_2,\ldots,Z_{n-1}$ the images of $x_1,x_2,\ldots,x_{n-1}$ under $\phi$,
 by $Z$ a general element of $A$, and consider the equation
\bes
f(Z_1,Z_2,\ldots,Z_{n-1},Z)=0
\ees
in $A$. Using the definition of the multiplication in $A$, we can rewrite the last equation in the form
\bee\label{f7}
h(x,\partial^{\a_1}(Z),\partial^{\a_2}(Z),\ldots,\partial^{\a_r}(Z))=0
\eee
where $h=h(x,t_{\a_1},\ldots,t_{\a_r})$ is a polynomial in the variables
$x,t_{\a_1},\ldots,t_{\a_r}$. Since $f\notin
\mathrm{N}\langle x_1,\ldots,x_{n-1}\rangle$ the polynomial $h$ essentially depends on $t_{\a_1},\ldots,t_{\a_r}$, i.e.
$r>0$ in (\ref{f4}).

Assume that $\a_1< \ldots <\a_r$ and that $h$ is irreducible. If $h$ is not irreducible we can replace it with its irreducible factor which contains $t_{\a_r}$.
We assert that there exists $L=(c,c_{\a_1},\ldots,c_{\a_r})\in k^{1+r}$
such that $h(L)=0$ and $\frac{\partial h}{\partial t_{\a_r}}(L)\neq 0$. If this is not true then by Hilbert's Nulstellenssatz $h$ divides $(\frac{\partial h}{\partial t_{\a_r}})^s$ for some $s>0$. But then, since $h$ is irreducible, $h$ divides $(\frac{\partial h}{\partial t_{\a_r}})$, which is clearly impossible.

Therefore we can use Corollary \ref{c6} and find a solution $Z_n$ of the differential equation (\ref{f7}) in the formal power series algebra $k[[x-c]]$. Note that $B=\langle k[[x-c]],\circ\rangle$ is a Novikov algebra and $A$ is a subalgebra of $B$.
Take a homomorphism of Novikov algebras  $\theta : \mathrm{N}\langle x_1,\ldots,x_n\rangle\rightarrow B$ defined by
\bes
\theta(x_1)=Z_1,\theta(z_2)=Z_2,\ldots,\theta(z_{n-1})=Z_{n-1},\theta(x_n)=Z_n.
\ees
Then $\theta_{|\mathrm{N}\langle x_1,\ldots,x_{n-1}\rangle}=\phi_{|\mathrm{N}\langle x_1,\ldots,x_{n-1}\rangle}$ and $\theta(f)=0$.
$\Box$

In many cases the Freiheitssatz is formulated directly in the language of freeness.
\begin{co}\label{c7}{\bf (Freiheitssatz)} Let $\mathrm{N}\langle x_1,\ldots,x_n\rangle$ be the free Novikov algebra over a field $k$ of characteristic $0$ in the variables $x_1,\ldots,x_n$. Suppose that $f \in \mathrm{N}\langle x_1,\ldots,x_n\rangle$ and $f\notin
\mathrm{N}\langle x_1,\ldots,x_{n-1}\rangle$. Then the subalgebra of the quotient algebra $\mathrm{N}\langle x_1,\ldots,x_n\rangle/\mathrm{id}(f)$ generated by $x_1+\mathrm{id}(f),\ldots,x_{n-1}+\mathrm{id}(f)$ is a free Novikov algebra with free generators $x_1+\mathrm{id}(f),\ldots,x_{n-1}+\mathrm{id}(f)$.
\end{co}

\bigskip

\begin{center}
{\bf\large Acknowledgments}
\end{center}

\hspace*{\parindent}

The authors are grateful to Professor Askar Dzhumadil'daev for interesting discussions.


\begin{thebibliography}{99}

\bibitem{Drensky00}
Drensky, V., Free algebras and PI-algebras,
Graduate course in algebra, Springer-Verlag Singapore, Singapore, 2000.

\bibitem{Novikov}   Balinskii, A.A.,   Novikov, S.P., Poisson bracket
of hamiltonian type, Frobenius algebras and Lie algebras. Dokladu
AN SSSR 283 (1985), No. 5, 1036--1039.

\bibitem{Dzhuma02}  Dzhumadil'daev,  A.S., Lofwall, C.,
Trees, free right-symmetric algebras, free Novikov algebras
and identities. Homology, Homotopy and Appl.  4 (2002),
No.2(1), 165-190.

\bibitem{Dzhuma05} Dzhumadil'daev, A.S., Tulenbaev, K.M., Nilpotency of Zinbiel algebras. J. Dyn. Control Syst. 11 (2005), no. 2, 195--213.

\bibitem{Dzhuma09} Dzhumadil'daev, A.S., Codimension growth and non-Koszulity of Novikov operad. arXiv:0902.3187, 7 pages.


\bibitem{Gelfand}   Gelfand, I.M.,  Dorfman, I.Ya., Hamiltonian operators and related algebraic structures. Funct. Anal. Appl. 13 (1979), 248--262.



\bibitem{KMLU}
Kozybaev, D., Makar-Limanov, L., Umirbaev, U., The Freiheitssatz and the automorphisms
 of free right-symmetric algebras. Asian-European Journal of Mathematics,  1 (2008), No. 2, 243--254.


\bibitem{Magnus} Magnus, M., \"Uber discontinuierliche Gruppen mit einer definierenden
Relation (Der Freiheitssatz). J. Reine Angew. Math. 163 (1930), 141–-165



\bibitem{Makar2} Makar-Limanov, L., Algebraically closed skew fields.  J. Algebra 93 (1985), no. 1, 117--135.



\bibitem{MLU2010} Makar-Limanov, L., Umirbaev, U., The Freiheitssatz
 for Poisson algebras. J. Algebra 328 (2011), 495-503.



\bibitem{NU}
Naurazbekova, A., Umirbaev, U., Identities of dual Leibniz algebras.
 TWMS J. Pure Appl. Math. 1 (2010), No. 1, 86--91.



\bibitem{Osborn}  Osborn, J.M., Infinite-dimensional Novikov algebras of
characteristic $0$. J. Algebra 167 (1994), 146--167.


\bibitem{Romanovskii}
Romanovskii, N.S., A theorem on freeness for groups with one defining relation in varieties of solvable and nilpotent groups of given degrees. (Russian) Mat. Sb. (N.S.) 89(131) (1972), 93--99.


\bibitem{Shir2}
Shirshov, A.I., Some algorithm problems for Lie algebras.
Sibirsk. Mat. Z. 3 (1962), 292--296.

\bibitem{KBKA} Zhevlakov, K.A., Slinko, A.M., Shestakov, I.P., Shirshov, A.I.,
 Rings that are nearly associative, Academic Press, Inc. [Harcourt Brace Jovanovich, Publishers], New York-London, 1982.

\end{thebibliography}
\end{document}